\documentclass[a4paper,11pt]{amsart}

\usepackage{graphicx}
\usepackage{mathptmx}
\usepackage{amsmath}
\usepackage{amssymb}
\usepackage{enumitem}
\usepackage{xcolor}

\newmuskip\pFqmuskip

\newcommand*\pFq[6][8]{%
  \begingroup 
  \pFqmuskip=#1mu\relax
  \mathcode`=\string"8000
  \begingroup\lccode`\~=`\,
  \lowercase{\endgroup\let~}\pFqcomma
  F^{#2}_{#3}{\left(\genfrac..{0pt}{}{#4}{#5}\bigg|#6\right)}%
  \endgroup
}
\newcommand{\pFqcomma}{\mskip\pFqmuskip}

\newtheorem{theorem}{Theorem}[section]

\newtheorem{remark}[theorem]{Remark}

\begin{document}

\title[]{Heterogeneous Stirling numbers and heterogeneous Bell polynomials}

\author{Taekyun  Kim}
\address{Department of Mathematics, Kwangwoon University, Seoul 139-701, Republic of Korea}
\email{tkkim@kw.ac.kr}
\author{Dae San  Kim}
\address{Department of Mathematics, Sogang University, Seoul 121-742, Republic of Korea}
\email{dskim@sogang.ac.kr}


\begin{abstract}
This paper introduces a novel generalization of Stirling and Lah numbers, termed ``heterogeneous Stirling numbers," which smoothly interpolate between these classical combinatorial sequences. Specifically, we define heterogeneous Stirling numbers of the second and first kinds, demonstrating their convergence to standard Stirling numbers as $\lambda \rightarrow 0$ and to (signed) Lah numbers as 
$\lambda \rightarrow 1$. We derive fundamental properties, including generating functions, explicit formulas, and recurrence relations. Furthermore, we extend these concepts to heterogeneous Bell polynomials, obtaining analogous results such as generating function, combinatorial identity and Dobinski-like formula. Finally, we introduce and analyse heterogeneous $r$-Stirling numbers of the second kind and their associated $r$-Bell polynomials.
\end{abstract}

\maketitle

\markboth{\centerline{\scriptsize Heterogeneous Stirling numbers and heterogeneous Bell polynomials}} 
{\centerline{\scriptsize Taekyun Kim and Dae San Kim}}

\section{Introduction} 
This paper introduces the heterogeneous Stirling numbers of the second kind, denoted by $H_{\lambda}(n,k)$ (defined in \eqref{13}), which generalize both the classical Stirling numbers of the second kind, ${n \brace k}$ (obtained as $\lambda \rightarrow 0$), and the Lah numbers, $L(n,k)$ (obtained as $\lambda \rightarrow 1$), (see \eqref{14}). We derive their generating function, explicit expressions, and recurrence relations. By inversion, we define the heterogeneous Stirling numbers of the first kind, $G_{\lambda}(n,k)$ (defined in \eqref{19}), which reduce to the classical Stirling numbers of the first kind $S_{1}(n,k)$ as $\lambda \rightarrow 0$ and to $(-1)^{n-k}L(n,k)$ as $\lambda \rightarrow 1$, (see \eqref{20}). We also establish their generating function and recurrence relation. As a natural polynomial extension of $H_{\lambda}(n,k)$, we introduce the heterogeneous Bell polynomials, $H_{n,\lambda}(x)$ (defined in \eqref{25}), and obtain their generating function, a combinatorial identity, and a Dobinski-like formula. Furthermore, we extend these concepts to the $r$-case, defining heterogeneous $r$-Stirling numbers of the second kind, $H_{\lambda}^{(r)}(n+r,k+r)$, (see \eqref{42}), and heterogeneous $r$-Bell polynomials, $H_{n,\lambda}^{(r)}(x)$, (see \eqref{47}), for which we derive analogous results. \par
The structure of this paper is as follows: Section 1 provides a review of degenerate exponentials, Stirling numbers of the first and second kinds, Lah numbers, Bell polynomials, Lah-Bell polynomials, degenerate Stirling numbers, degenerate Bell polynomials, and degenerate $r$-Stirling numbers. Sections 2 and 3 present the main contributions of this work. Specifically, for $H_{\lambda}(n,k)$, we derive a generating function (Theorem 2.1), explicit expressions (Theorems 2.2 and 2.10), and a recurrence relation (Theorem 2.5). For $G_{\lambda}(n,k)$, we establish a generating function (Theorem 2.4) and a recurrence relation (Theorem 2.7). For $H_{n,\lambda}(x)$, we obtain a generating function (Theorem 2.8), a Dobinski-like formula (Theorem 2.12), and a combinatorial identity (Theorem 3.1). Additionally, we present generating functions and recurrence relations for $H_{\lambda}^{(r)}(n+r,k+r)$ (\eqref{45} and \eqref{46}, respectively) and generating functions and explicit expressions for $H_{n,\lambda}^{(r)}(x)$ (\eqref{48} and \eqref{50}, respectively). One may see [2,6,7,18,19] as general references. The remainder of this section provides essential background information used throughout the paper. \par
\vspace{0.1in}
The falling factorial sequence is defined by (see [4,5,9-15])
\begin{equation*}
(x)_{0}=1,\quad (x)_{n}=x(x-1)\cdots(x-n+1),\quad (n\ge 1),
\end{equation*}
and the rising factorial sequence is given by 
\begin{equation*}
\langle x\rangle_{0}=1,\quad \langle x\rangle_{n}=x(x+1)\cdots(x+n-1),\quad (n\ge 1).
\end{equation*}
For any nonzero $\lambda\in\mathbb{R}$, the degenerate falling factorial sequence is defined by (see [16])
\begin{equation*}
(x)_{0,\lambda}=1,\quad (x)_{n,\lambda}=x(x-\lambda)(x-2\lambda)\cdots\big(x-(n-1)\lambda\big),\quad (n\ge 1), 
\end{equation*}
and the degenerate rising factorial sequence is given by
\begin{equation*}
\langle x\rangle_{0,\lambda}=1,\quad \langle x\rangle_{n,\lambda}=x(x+\lambda)(x+2\lambda)\cdots\big(x+(n-1)\lambda\big),\quad (n\ge 1).
\end{equation*}
Note that 
\begin{displaymath}
\lim_{\lambda\rightarrow 0}(x)_{n,\lambda}=\langle x\rangle_{n,\lambda}=x^{n},\ \lim_{\lambda\rightarrow 1}(x)_{n,\lambda}=(x)_{n},\ \lim_{\lambda\rightarrow 1}\langle x\rangle_{n,\lambda}=\langle x\rangle_{n}.
\end{displaymath} \par
The degenerate exponentials are defined by 
\begin{equation}
e_{\lambda}^{x}(t)=\sum_{n=0}^{\infty}(x)_{n,\lambda}\frac{t^{n}}{n!},\quad e_{\lambda}(t)=e_{\lambda}^{1}(t),\quad (\mathrm{see}\ [9-16]). \label{1}
\end{equation}
The Stirling numbers of the first kind are given by 
\begin{equation*}
(x)_{n}=\sum_{k=0}^{n}S_{1}(n,k)x^{k},\quad (n\ge 0),\quad (\mathrm{see}\ [3,4,6,8,18,19]).
\end{equation*}
The unsigned Stirling numbers of the first kind are defined by 
\begin{align*}
&{n \brack k}=(-1)^{n-k}S_{1}(n,k),\quad (n \ge k\ge 0),\quad (\mathrm{see}\ [6,16,19]),  \\
&\frac{1}{k!}\big(-\log (1-t) \big)^{k}=\sum_{n=k}^{\infty} {n \brack k} \frac{t^{n}}{n!},\quad (k \ge 0). 
\end{align*}
The Stirling numbers of the second kind are given by   
\begin{equation*}
x^{n}=\sum_{k=0}^{n}{n \brace k}(x)_{k},\quad (n\ge 0),\quad (\mathrm{see}\ [3,4,6,8,18,19]). 
\end{equation*} \par
The Lah numbers are defined by 
\begin{equation}
\frac{1}{k!}\bigg(\frac{t}{1-t}\bigg)^{k}=\sum_{n=k}^{\infty}L(n,k)\frac{t^{n}}{n!},\ (k\ge 0),\ (\mathrm{see}\ [6,11,12,19]). \label{2}	
\end{equation}
Thus, by \eqref{2}, we get 
\begin{equation}
L(n,k)=\frac{n!}{k!}\binom{n-1}{k-1},\quad (n \ge k\ge 0).\label{3}
\end{equation}
The Bell polynomials are given by 
\begin{equation}
\phi_{n}(x)=\sum_{k=0}^{n}{n \brace k}x^{k},\quad (n\ge 0),\quad (\mathrm{see}\ [6,13,14,16,19]).\label{4}
\end{equation}
When $x=1$, $\phi_{n}=\phi_{n}(1)$ are called the Bell numbers, (see [1,20]). \\
Thus, by \eqref{4}, we get 
\begin{equation}
e^{x(e^{t}-1)}=\sum_{n=0}^{\infty}\phi_{n}(x)\frac{t^{n}}{n!},\quad (\mathrm{see}\ [6,13,14,16,19]).\label{5}
\end{equation}
From \eqref{5}, we obtain the following Dobinski's formula:
\begin{equation}
\phi_{n}(x)=e^{-x}\sum_{k=0}^{\infty}\frac{k^{n}}{k!}x^{k}. \label{6}
\end{equation} \\
The Lah-Bell polynomials are defined by 
\begin{equation}
e^{x\big(\frac{1}{1-t}-1\big)}=\sum_{n=0}^{\infty}\mathrm{LB}_{n}(x)\frac{t^{n}}{n!},\quad (\mathrm{see}\ [6,11,12,17,19]).\label{7}
\end{equation}
When $x=1$, $LB_{n}=LB_{n}(1)$ are called the Lah-Bell numbers. \\
From \eqref{7}, we derive the next Dobinski-like formula:
\begin{equation}
LB_{n}(x)=e^{-x}\sum_{k=0}^{\infty}\frac{\langle k \rangle_{n}}{k!}x^{k}. \label{8}
\end{equation}
Thus, by \eqref{2} and \eqref{7}, we get 
\begin{equation*}
\mathrm{LB}_{n}(x)=\sum_{k=0}^{n}x^{k}L(n,k),\quad (n\ge 0),\quad (\mathrm{see}\ [11,12]).
\end{equation*} \par
Let $\log_{\lambda}t$ be the compositional inverse function of $e_{\lambda}(t)$. Then we have 
\begin{equation*}
\log_{\lambda}(1+t)=\sum_{n=1}^{\infty}(1)_{n,1/\lambda}\lambda^{n-1}\frac{t^{n}}{n!},\quad (\mathrm{see}\ [9]),
\end{equation*}
and 
\begin{equation*}
e_{\lambda}\log_{\lambda} t=\log_{\lambda}e_{\lambda}(t)=t.
\end{equation*}
The unsigned degenerate Stirling numbers of the first kind are given by 
\begin{align}
&\langle x\rangle_{n}=\sum_{k=0}^{n}{n \brack k}_{\lambda}\langle x\rangle_{k,\lambda},\quad (n\ge 0),\quad (\mathrm{see}\ [16]), \label{9}	\\
&\frac{1}{k!}\big(-\log_{\lambda} (1-t) \big)^{k}=\sum_{n=k}^{\infty} {n \brack k}_{\lambda} \frac{t^{n}}{n!},\quad (k \ge 0). \nonumber
\end{align}
The degenerate Stirling numbers of the first kind are defined by 
\begin{displaymath}
S_{1,\lambda}(n,k)=(-1)^{n-k}{n \brack k}_{\lambda},\quad (n \ge k\ge 0). 
\end{displaymath}
Inversion of \eqref{9} gives the degenerate Stirling numbers of the second kind defined by 
\begin{align}
&(x)_{n,\lambda}=\sum_{k=0}^{n}{n \brace k}_{\lambda}(x)_{k},\quad (n\ge 0),\quad (\mathrm{see}\ [9,14-17]),\label{10} \\
&\frac{1}{k!}\big(e_{\lambda}(t)-1 \big)^{k}=\sum_{n=k}^{\infty} {n \brace k}_{\lambda}\frac{t^{n}}{n!}, \quad (k \ge 0). \nonumber
\end{align}
The degenerate Bell polynomials are defined by 
\begin{equation}
\phi_{n,\lambda}(x)=\sum_{k=0}^{n}{n \brace k}_{\lambda}x^{k},\quad (n\ge 0),\quad (\mathrm{see}\ [14,15,16]).\label{11}
\end{equation}
Thus, by \eqref{10} and \eqref{11}, we get 
\begin{equation*}
e^{x(e_{\lambda}(t)-1)}=\sum_{n=0}^{\infty}\phi_{n,\lambda}(x)\frac{t^{n}}{n!}. 
\end{equation*} \par
For $r\in\mathbb{N}$, the degenerate $r$-Stirling numbers of the second kind are given by 
\begin{equation}
(x+r)_{n,\lambda}=\sum_{k=0}^{n}{n+r \brace k+r}_{r,\lambda}(x)_{k},\quad (n\ge 0),\quad (\mathrm{see}\ [16]).\label{12}
\end{equation}
Inversion of \eqref{12} yields the unsigned degenerate $r$-Stirling numbers of the first kind defined by 
\begin{equation*}
\langle x+r\rangle_{n}=\sum_{k=0}^{n}{n+r \brack k+r}_{r,\lambda}\langle x\rangle_{k,\lambda},\quad (n\ge 0),\quad (\mathrm{see}\ [16]). 
\end{equation*}

\section{Heterogeneous Stirling numbers and heterogeneous Bell polynomials}
In this section, we introduce and investigate the heterogeneous Stirling numbers of the second kind $H_{\lambda}(n,k)$ and of the first kind $G_{\lambda}(n,k)$, along with a polynomial extension of $H_{\lambda}(n,k)$, the heterogeneous Bell polynomials $H_{n,\lambda}(x)$. \par
We define the {\it{heterogeneous Stirling numbers of the second kind}} by
\begin{equation}
\langle x\rangle_{n,\lambda}=\sum_{k=0}^{n}H_{\lambda}(n,k)(x)_{k},\quad (n\ge 0). \label{13}
\end{equation}
As usual, we define $H_{\lambda}(n,k)=0$,\,\, for $k <0$\, or\, $n < k$.
Note that 
\begin{equation}
\lim_{\lambda\rightarrow 0}H_{\lambda}(n,k)={n \brace k},\quad\mathrm{and}\quad \lim_{\lambda\rightarrow 1}H_{\lambda}(n,k)=L(n,k),\quad (n \ge k \ge 0). \label{14}
\end{equation}
From \eqref{13}, we note that 
\begin{align}
e_{\lambda}^{-x}(-t)&=\sum_{n=0}^{\infty}\frac{1}{n!}\langle x\rangle_{n,\lambda}t^{n}=\sum_{n=0}^{\infty}\frac{t^{n}}{n!}\sum_{k=0}^{n}H_{\lambda}(n,k)(x)_{k}\label{15}\\
&=\sum_{k=0}^{\infty}\sum_{n=k}^{\infty}H_{\lambda}(n,k)\frac{t^{n}}{n!}(x)_{k}. \nonumber
\end{align}
On the other hand, by binomial theorem, we get 
\begin{equation}
e_{\lambda}^{-x}(-t)=\big(e_{\lambda}^{-1}(-t)-1+1\big)^{x}=\sum_{k=0}^{\infty}\frac{1}{k!}\big(e_{\lambda}^{-1}(-t)-1\big)^{k}(x)_{k}. \label{16}
\end{equation}
Therefore, by \eqref{15} and \eqref{16}, we obtain the following theorem. 
\begin{theorem}
For $k\ge 0$, the generating function of heterogeneous Stirling numbers of the second kind is given by
\begin{equation*}
\frac{1}{k!}\big(e_{\lambda}^{-1}(-t)-1\big)^{k}=\sum_{n=k}^{\infty}H_{\lambda}(n,k)\frac{t^{n}}{n!}.
\end{equation*}
\end{theorem}
From \eqref{15} and Theorem 2.1, we note that 
\begin{align}
\frac{1}{k!}\big(e_{\lambda}^{-1}(-t)-1\big)^{k}&=\frac{1}{k!}\sum_{j=0}^{k}\binom{k}{j}(-1)^{k-j}e_{\lambda}^{-j}(-t)\label{17}\\
&=\sum_{n=0}^{\infty}\frac{1}{k!}\sum_{j=0}^{k}\binom{k}{j}(-1)^{k-j}\langle j\rangle_{n,\lambda}\frac{t^{n}}{n!}.\nonumber	
\end{align}
Thus, by Theorem 2.1 and \eqref{17}, we get 
\begin{equation}
\frac{1}{k!}\sum_{j=0}^{k}\binom{k}{j}(-1)^{k-j}\langle j\rangle_{n,\lambda}=\left\{\begin{array}{ccc}
H_{\lambda}(n,k), & \textrm{if $n\ge k$,} \\
0, & \textrm{if $0 \le n<k$.}
\end{array}\right.\label{18}
\end{equation}
Therefore, by \eqref{18}, we obtain the following theorem.
\begin{theorem}
For $n\ge k\ge 0$, we have
\begin{displaymath}
H_{\lambda}(n,k)=\frac{1}{k!}\sum_{j=0}^{k}\binom{k}{j}(-1)^{k-j}\langle j\rangle_{n,\lambda}.
\end{displaymath}
\end{theorem}
\begin{remark}
In Theorem 2.2, letting $\lambda \rightarrow 1$, and letting $\lambda \rightarrow 0$ respectively yield
\begin{align*}
&L(n,k)=\frac{1}{k!}\sum_{j=0}^{k}\binom{k}{j}(-1)^{k-j}\langle j\rangle_{n}, \\
&{n \brace k}=\frac{1}{k!}\sum_{j=0}^{k}\binom{k}{j}(-1)^{k-j}j^{n}. 
\end{align*}
\end{remark}
Now, we consider the {\it{heterogeneous Stirling numbers of the first kind}} which are defined by 
\begin{equation}
(x)_{n}=\sum_{k=0}^{n}G_{\lambda}(n,k)\langle x\rangle_{k,\lambda},\quad (n\ge 0). \label{19}
\end{equation}
Note that 
\begin{equation}
\lim_{\lambda\rightarrow 0}G_{\lambda}(n,k)=S_{1}(n,k)\quad\mathrm{and}\quad \lim_{\lambda\rightarrow 1}G_{\lambda}(n,k)=(-1)^{n-k}L(n,k),\quad (n \ge k\ge 0).  \label{20}
\end{equation}
From \eqref{19}, we have 
\begin{align}
(1+t)^{x}&=\sum_{n=0}^{\infty}(x)_{n}\frac{t^{n}}{n!}=\sum_{n=0}^{\infty}\sum_{k=0}^{n}G_{\lambda}(n,k)\langle x\rangle_{k,\lambda}\frac{t^{n}}{n!}\label{21}\\
&=\sum_{k=0}^{\infty}\sum_{n=k}^{\infty}G_{\lambda}(n,k)\frac{t^{n}}{n!}\langle x\rangle_{k,\lambda}.\nonumber
\end{align}
On the other hand, by binomial expansion, we get 
\begin{align}
(1+t)^{x}&=e_{-\lambda}^{x}\Big(\log_{-\lambda}(1+t)\Big)=\sum_{k=0}^{\infty}\frac{(x)_{k,-\lambda}}{k!}\Big(\log_{-\lambda}(1+t)\Big)^{k} \label{22}\\
&=\sum_{k=0}^{\infty}\frac{1}{k!}\Big(\log_{-\lambda}(1+t)\Big)^{k}\langle x\rangle_{k,\lambda}.\nonumber	
\end{align}
Therefore, by \eqref{21} and \eqref{22}, we obtain the following theorem. 
\begin{theorem}
For $k\ge 0$, the generating function of heterogeneous Stirling numbers of the first kind is given by
\begin{displaymath}
\frac{1}{k!}\Big(\log_{-\lambda}(1+t)\Big)^{k}=\sum_{n=k}^{\infty}G_{\lambda}(n,k)\frac{t^{n}}{n!}. 
\end{displaymath}
\end{theorem}
From \eqref{13}, we note that 
\begin{align}
&\sum_{k=0}^{n+1}H_{\lambda}(n+1,k)(x)_{k}=\langle x\rangle_{n+1,\lambda}=\langle x\rangle_{n,\lambda}(x+n\lambda)\label{23}\\
&=\sum_{k=0}^{n}H_{\lambda}(n,k)(x)_{k}(x-k+k+n\lambda)\nonumber\\
&=\sum_{k=0}^{n}H_{\lambda}(n,k)(x)_{k+1}+\sum_{k=0}^{n}H_{\lambda}(n,k)(k+n\lambda)(x)_{k}\nonumber\\
&=\sum_{k=1}^{n+1}H_{\lambda}(n,k-1)(x)_{k}+\sum_{k=0}^{n}H_{\lambda}(n,k)(k+n\lambda)(x)_{k}\nonumber\\
&=\sum_{k=0}^{n+1}\Big(H_{\lambda}(n,k-1)+(k+n\lambda)H_{\lambda}(n,k)\Big)(x)_{k}.\nonumber
\end{align}
Therefore, by comparing coefficients on both sides of \eqref{23}, we obtain the following theorem. 
\begin{theorem}
For $n \ge k\ge 0$, we have 
\begin{displaymath}
H_{\lambda}(n+1,k)=H_{\lambda}(n,k-1)+(k+n\lambda)H_{\lambda}(n,k).
\end{displaymath}
\end{theorem}
\begin{remark} 
In Theorem 2.5, letting $\lambda \rightarrow 1$, and letting $\lambda \rightarrow 0$ respectively give
\begin{align*}
&L(n+1,k)=L(n,k-1)+(k+n)L(n,k), \\
&{n+1 \brace k}={n \brace k-1}+k{n \brace k}. 
\end{align*}
\end{remark}
By \eqref{19}, we get 
\begin{align}
&\sum_{k=0}^{n+1}G_{\lambda}(n+1,k)\langle x\rangle_{k,\lambda}=(x)_{n+1}=(x)_{n}(x-n) \label{24}\\
&=\sum_{k=0}^{n}G_{\lambda}(n,k)\langle x\rangle_{k,\lambda}(x+k\lambda-n-k\lambda)\nonumber \\
&=\sum_{k=0}^{n}G_{\lambda}(n,k)\langle x\rangle_{k+1,\lambda}+\sum_{k=0}^{n}(-n-k\lambda)G_{\lambda}(n,k)\langle x\rangle_{k,\lambda}\nonumber\\
&=\sum_{k=1}^{n+1}G_{\lambda}(n,k-1)\langle x\rangle_{k,\lambda}-\sum_{k=0}^{n}(n+k\lambda)G_{\lambda}(n,k)\langle x\rangle_{k,\lambda} \nonumber \\
&=\sum_{k=0}^{n+1}\Big(G_{\lambda}(n,k-1)-(n+k\lambda)G_{\lambda}(n,k)\Big)\langle x\rangle_{k,\lambda}.\nonumber	
\end{align}
Therefore, by comparing the coefficients on both sides of \eqref{24}, we obtain the following theorem. 
\begin{theorem}
For $n \ge k\ge 0$, we have 
\begin{equation*}
G_{\lambda}(n+1,k)=G_{\lambda}(n,k-1)-(n+k\lambda)G_{\lambda}(n,k).
\end{equation*}
\end{theorem}
\begin{remark}
In Theorem 2.7, letting $\lambda \rightarrow 1$, and letting $\lambda \rightarrow 0$ respectively give
\begin{align*}
&L(n+1,k)=L(n,k-1)+(n+k)L(n,k), \\
&S_{1}(n+1,k)=S_{1}(n,k-1)-nS_{1}(n,k).
\end{align*}
\end{remark}
In view of \eqref{11}, we define the {\it{heterogeneous Bell polynomials}} by
\begin{equation}
H_{n,\lambda}(x)=\sum_{k=0}^{n}H_{\lambda}(n,k)x^{k},\quad (n\ge 0). \label{25}	
\end{equation}
When $x=1$, $H_{n,\lambda}=H_{n,\lambda}(1)$ are called the {\it{heterogeneous Bell numbers}}. \par 
By Theorem 2.1 and \eqref{25}, we get 
\begin{align}
&\sum_{n=0}^{\infty}H_{n,\lambda}(x)\frac{t^{n}}{n!}=\sum_{n=0}^{\infty}\sum_{k=0}^{n}H_{\lambda}(n,k)x^{k}\frac{t^{n}}{n!}\label{26} \\
&=\sum_{k=0}^{\infty}\sum_{n=k}^{\infty}H_{\lambda}(n,k)\frac{t^{n}}{n!}x{^k}=\sum_{k=0}^{\infty}\frac{1}{k!}\Big(e_{\lambda}^{-1}(-t)-1\Big)^{k}x{^k}\nonumber\\
&=e^{x(e_{\lambda}^{-1}(-t)-1)}.\nonumber
\end{align}
Therefore, by \eqref{26}, we obtain the following theorem. 
\begin{theorem}
The generating function of the heterogeneous Bell polynomials is given by
\begin{displaymath}
e^{x(e_{\lambda}^{-1}(-t)-1)}= \sum_{n=0}^{\infty}H_{n,\lambda}(x)\frac{t^{n}}{n!}.
\end{displaymath}
\end{theorem}
From Theorem 2.1, \eqref{1} and \eqref{2}, we note that 
\begin{align}
&\sum_{n=k}^{\infty}H_{\lambda}(n,k)\frac{t^{n}}{n!}=\frac{1}{k!}\Big(e_{\lambda}^{-1}(-t)-1\Big)^{k}=\frac{1}{k!}\Big((1-\lambda t)^{-1/\lambda}-1\Big)^{k} \label{27} \\
&=\frac{1}{k!}\Big(e^{-\frac{1}{\lambda}\log(1-\lambda t)}-1\Big)^{k}=\sum_{l=k}^{\infty}{l \brace k}\lambda^{-l}\frac{1}{l!}\Big(-\log(1-\lambda t)\Big)^{l}\nonumber\\
&=\sum_{l=k}^{\infty}{l \brace k}\lambda^{-l}\sum_{n=l}^{\infty}{n \brack l}\lambda^{n}\frac{t^{n}}{n!}=\sum_{n=k}^{\infty}\sum_{l=k}^{n}{l \brace k}{n \brack l}\lambda^{n-l}\frac{t^{n}}{n!}.\nonumber
\end{align}
Therefore, by \eqref{27}, we obtain the following theorem. 
\begin{theorem}
For $n \ge k\ge 0$, we have 
\begin{displaymath}
H_{\lambda}(n,k)=\sum_{l=k}^{n}{l \brace k}{n \brack l}\lambda^{n-l}. 
\end{displaymath}
\end{theorem}
\begin{remark}
In Theorem 2.10, letting $\lambda \rightarrow 1$ yields
\begin{equation}
L(n,k)=\sum_{l=k}^{n}{l \brace k}{n \brack l}. \label{28}
\end{equation}
\end{remark}
From \eqref{10}, we note that 
\begin{equation}
\frac{1}{k!}\big(e_{\lambda}(t)-1\big)^{k}=\sum_{n=k}^{\infty}{n \brace k}_{\lambda}\frac{t^{n}}{n!},\quad (k\ge 0).\label{29}
\end{equation}
Replacing $t$ by $\log_{\lambda}\big(\frac{1}{1-t}\big)$ in \eqref{29} and using \eqref{2}, we have 
\begin{align}
&\sum_{n=k}^{\infty}L(n,k)\frac{t^{n}}{n!}=\frac{1}{k!}\bigg(\frac{1}{1-t}-1\bigg)^{k}=\sum_{l=k}^{\infty}{l \brace k}_{\lambda}\frac{1}{l!}\bigg(\log_{\lambda}\bigg(\frac{1}{1-t}\bigg)\bigg)^{l} \label{30}\\
&=\sum_{l=k}^{\infty}{l \brace k}_{\lambda}\frac{1}{l!}\Big(-\log_{-\lambda}(1-t)\Big)^{l}=\sum_{l=k}^{\infty}{l \brace k}_{\lambda}\sum_{n=l}^{\infty}{n \brack l}_{-\lambda}\frac{t^{n}}{n!} \nonumber\\
&=\sum_{n=k}^{\infty}\sum_{l=k}^{n}{l \brace k}_{\lambda}{n \brack l}_{-\lambda}\frac{t^{n}}{n!}.\nonumber
\end{align}
Therefore, by comparing the coefficients on both sides of \eqref{30}, we obtain the following theorem. 
\begin{theorem}
For $n\ge k\ge 0$, we have 
\begin{displaymath}
L(n,k)=\sum_{l=k}^{n}{l \brace k}_{\lambda}{n \brack l}_{-\lambda}.
\end{displaymath}
Equivalently, we have
\begin{displaymath}
\binom{n}{k}\binom{n-1}{k-1}=\frac{1}{(n-k)!}\sum_{l=k}^{n}{l \brace k}_{\lambda}{n \brack l}_{-\lambda}, \quad (\textnormal{see}\,\,\eqref{3}).
\end{displaymath}
\end{theorem}
Replacing $t$ by $-\log_{\lambda}(1-t)$ in Theorem 2.1 and using \eqref{2} and \eqref{9}, we get 
\begin{align}
&\sum_{n=k}^{\infty}L(n,k)\frac{t^{n}}{n!}=\frac{1}{k!}\bigg(\frac{1}{1-t}-1\bigg)^{k}=\sum_{l=k}^{\infty}H_{\lambda}(l,k)\frac{1}{l!}\Big(-\log_{\lambda}(1-t)\Big)^{l}\label{31}\\
&=\sum_{l=k}^{\infty}H_{\lambda}(l,k)\sum_{n=l}^{\infty}{n \brack l}_{\lambda}\frac{t^{n}}{n!}=\sum_{n=k}^{\infty}\sum_{l=k}^{n}H_{\lambda}(l,k){n \brack l}_{\lambda}\frac{t^{n}}{n!}.\nonumber
\end{align}
Therefore, by \eqref{31}, we obtain the following theorem. 
\begin{theorem}
For $n \ge k\ge 0$, we have 
\begin{displaymath}
L(n,k)=\sum_{l=k}^{n}H_{\lambda}(l,k){n \brack l}_{\lambda}. 
\end{displaymath}
Equivalently, we have
\begin{displaymath}
\binom{n}{k}\binom{n-1}{k-1}=\frac{1}{(n-k)!}\sum_{l=k}^{n}H_{\lambda}(l,k){n \brack l}_{\lambda}\quad (\textnormal{see}\,\,\eqref{3}). 
\end{displaymath}
\end{theorem}
\begin{remark}
Letting $\lambda \rightarrow 0$ in Theorem 2.13 yields the same identity as in \eqref{28}:
\begin{equation*}
L(n,k)=\sum_{l=k}^{n}{l \brace k}{n \brack l}.
\end{equation*}
\end{remark}
From Theorem 2.9, we note that 
\begin{align}
\sum_{n=0}^{\infty}H_{n,\lambda}(x)\frac{t^{n}}{n!}&=e^{-x}e^{xe_{\lambda}^{-1}(-t)}=e^{-x}\sum_{k=0}^{\infty}\frac{x^{k}}{k!}e_{\lambda}^{-k}(-t) \label{32}\\
&=\sum_{n=0}^{\infty}e^{-x}\sum_{k=0}^{\infty}\frac{x^{k}}{k!}\langle k\rangle_{n,\lambda}\frac{t^{n}}{n!}. \nonumber	
\end{align}
Therefore, by \eqref{32}, we obtain the following Dobinski-like theorem. 
\begin{theorem}
For $n\ge 0$, we have 
\begin{displaymath}
H_{n,\lambda}(x)=e^{-x}\sum_{k=0}^{\infty}\frac{\langle k\rangle_{n,\lambda}}{k!}x^{k}. 
\end{displaymath}
\end{theorem}
\begin{remark}
In Theorem 2.15, letting $\lambda \rightarrow 1$, and letting $\lambda \rightarrow 0$ respectively give
\begin{equation*}
\mathrm{LB}_{n}(x)=e^{-x}\sum_{k=0}^{\infty}\frac{\langle k\rangle_{n}}{k!}x^{k}, \quad
\phi_{n}(x)=e^{-x}\sum_{k=0}^{\infty}\frac{k^{n}}{k!}x^{k},\quad (\textnormal{see}\,\, \eqref{6},\, \eqref{8}). 
\end{equation*}
\end{remark}
From Theorem 2.9, we note that 
\begin{align}
&\sum_{m,n=0}^{\infty}H_{n+m,\lambda}\frac{x^{n}}{n!}\frac{y^{m}}{m!}=\sum_{m=0}^{\infty}\bigg(\frac{d}{dx}\bigg)^{m}\sum_{n=0}^{\infty}H_{n,\lambda}\frac{x^{n}}{n!}\frac{y^{m}}{m!}\label{33}	\\
&=\sum_{m=0}^{\infty}\frac{y^{m}}{m!}\bigg(\frac{d}{dx}\bigg)^{m}\sum_{n=0}^{\infty}H_{n,\lambda}\frac{x^{n}}{n!}=e^{y\frac{d}{dx}}e^{e_{\lambda}^{-1}(-x)-1}. \nonumber
\end{align}
Now, we observe that 
\begin{align}
&\sum_{m=0}^{\infty}\bigg(\frac{d}{dx}\bigg)^{m}\sum_{n=0}^{\infty}H_{n,\lambda}\frac{x^{n}}{n!}\frac{y^{m}}{m!}=\sum_{n=0}^{\infty}\frac{H_{n,\lambda}}{n!}\sum_{m=0}^{\infty}\bigg(\frac{d}{dx}\bigg)^{m}x^{n}\frac{y^{m}}{m!} \label{34}\\
&=\sum_{n=0}^{\infty}\frac{H_{n,\lambda}}{n!}(x+y)^{n}=e^{e_{\lambda}^{-1}(-x-y)-1}.\nonumber
\end{align}
Thus, by \eqref{33} and \eqref{34}, we get 
\begin{equation}
\sum_{m,n=0}^{\infty}H_{n+m,\lambda}\frac{x^{n}}{n!}\frac{y^{m}}{m!}=e^{y\frac{d}{dx}} e^{e_{\lambda}^{-1}(-x)-1}=e^{e_{\lambda}^{-1}(-x-y)-1}.\label{35}
\end{equation}
By \eqref{1}, we get 
\begin{align}
e_{\lambda}^{-1}(-x-y)&=\Big(1-\lambda(x+y)\Big)^{-\frac{1}{\lambda}}=\big(1-\lambda x-\lambda y\big)^{-\frac{1}{\lambda}}\label{36} \\
&=\big(1-\lambda x\big)^{-\frac{1}{\lambda}}\bigg(1-\frac{\lambda y}{1-\lambda x}\bigg)^{-\frac{1}{\lambda}}=e_{\lambda}^{-1}(-x)e_{\lambda}^{-1}\bigg(\frac{-y}{1-\lambda x}\bigg). \nonumber	
\end{align}
From \eqref{35} and \eqref{36} and using Theorems 2.1 and 2.9, we have 
\begin{align}
&\sum_{m,n=0}^{\infty}H_{m+n,\lambda}\frac{x^{n}}{n!}\frac{y^{m}}{m!}= e^{e_{\lambda}^{-1}(-x)-1} e^{e_{\lambda}^{-1}(-x)\big(e_{\lambda}^{-1}\big(\frac{-y}{1-\lambda x}\big)-1\big)}\label{37}\\
&= e^{e_{\lambda}^{-1}(-x)-1}\sum_{j=0}^{\infty}\frac{1}{j!}\bigg(e_{\lambda}^{-1}\bigg(\frac{-y}{1-\lambda x}\bigg)-1\bigg)^{j}e_{\lambda}^{-j}(-x)\nonumber\\
&= e^{e_{\lambda}^{-1}(-x)-1}\sum_{j=0}^{\infty}\sum_{m=j}^{\infty}H_{\lambda}(m,j)\frac{1}{m!}\bigg(\frac{y}{1-\lambda x}\bigg)^{m}e_{\lambda}^{-j}(-x)\nonumber\\
&= e^{e_{\lambda}^{-1}(-x)-1} \sum_{m=0}^{\infty}\sum_{j=0}^{m}H_{\lambda}(m,j)\frac{y^{m}}{m!}e_{\lambda}^{-(m\lambda+j)}(-x)\nonumber\\
&=\sum_{k=0}^{\infty}H_{k,\lambda}\frac{x^{k}}{k!}\sum_{m=0}^{\infty}\sum_{j=0}^{m}H_{\lambda}(m,j)\frac{y^{m}}{m!}e_{\lambda}^{-(m\lambda +j)}(-x)\nonumber\\
&=\sum_{k=0}^{\infty}H_{k,\lambda}\frac{x^{k}}{k!}\sum_{m=0}^{\infty}\sum_{j=0}^{m}H_{\lambda}(m,j)\frac{y^{m}}{m!}\sum_{l=0}^{\infty}\langle m\lambda +j\rangle_{l,\lambda}\frac{x^{l}}{l!}\nonumber\\
&=\sum_{m=0}^{\infty}\sum_{j=0}^{m}H_{\lambda}(m,j)\frac{y^{m}}{m!}\sum_{n=0}^{\infty}\sum_{k=0}^{n}\binom{n}{k}H_{k,\lambda}\langle m\lambda+ j\rangle_{n-k,\lambda}\frac{x^{n}}{n!} \nonumber\\
&=\sum_{m,n=0}^{\infty}\bigg(\sum_{j=0}^{m}H_{\lambda}(m,j)\sum_{k=0}^{n}\binom{n}{k}H_{k,\lambda}\langle m\lambda+j\rangle_{n-k,\lambda}\bigg)\frac{y^{m}}{m!}\frac{x^{n}}{n!}. \nonumber
\end{align}
By comparing the coefficients on both sides of \eqref{37}, we obtain the following theorem. 
\begin{theorem}
For $m,n\ge 0$, we have 
\begin{displaymath}
H_{m+n,\lambda}=\sum_{j=0}^{m}\sum_{k=0}^{n}\binom{n}{k}H_{\lambda}(m,j)H_{k,\lambda}\langle m\lambda+j\rangle_{n-k,\lambda}. 
\end{displaymath}
\end{theorem}
\begin{remark}
In Theorem 2.17, letting $ \lambda \rightarrow 0$, and letting $\lambda \rightarrow 1$ respectively give
\begin{align*}
&\phi_{m+n}=\sum_{j=0}^{m}\sum_{k=0}^{n}{m \brace j}\binom{n}{k}\phi_{k}j^{n-k},   \\
&\mathrm{LB}_{n+m}=\sum_{j=0}^{m}\sum_{k=0}^{n}L(m,j)\binom{n}{k}\mathrm{LB}_{k}\langle m+j\rangle_{n-k},
\end{align*}
where the first one is due to Spivey (see [21]).
\end{remark}

\section{Further remarks} 
In this section, we study the heterogeneous $r$-Stirling numbers of the second kind, $H_{\lambda}^{(r)}(n+r,k+r)$, and the heterogeneous $r$-Bell polynomials, $H_{n,\lambda}^{(r)}(x)$, respectively as $r$-analogues of $H_{\lambda}(n,k)$ and $H_{n,\lambda}(x)$. In this section, we use the following notations for operators:
\begin{equation*}
D=\frac{d}{dx}, \quad xD=x\frac{d}{dx}.
\end{equation*}
From Theorem 2.15, we note that 
\begin{align}
\langle xD\rangle_{n,\lambda}e^{x}&=\sum_{k=0}^{\infty}\frac{1}{k!}\langle xD\rangle_{n,\lambda}x^{k}=\sum_{k=0}^{\infty}\frac{\langle k\rangle_{n,\lambda}}{k!}x^{k} \label{38}\\
&=e^{x}\bigg(e^{-x}\sum_{k=0}^{\infty}\frac{\langle k\rangle_{n,\lambda}}{k!}x^{k}\bigg)=e^{x}H_{n,\lambda}(x).\nonumber
\end{align}
By \eqref{38}, we get 
\begin{align}
&H_{n+m,\lambda}(x)e^{x}=\langle xD\rangle_{m+n,\lambda}e^{x}=\langle xD+m\lambda\rangle_{n,\lambda}\langle xD\rangle_{m,\lambda}e^{x} \label{39}\\
&=\langle xD+m\lambda\rangle_{n,\lambda}\Big(H_{m,\lambda}(x)e^{x}\Big)=\sum_{k=0}^{n}\binom{n}{k}\Big[\langle xD \rangle_{k,\lambda}e^{x} \Big]\Big[\langle xD+m \lambda \rangle_{n-k,\lambda}H_{m,\lambda}(x) \Big]\nonumber\\
&=\sum_{k=0}^{n}\binom{n}{k}H_{k,\lambda}(x)\sum_{j=0}^{m}H_{\lambda}(m,j)\langle j+m\lambda \rangle_{n-k,\lambda}x^{j} e^{x} \nonumber \\
&=\sum_{k=0}^{n}\sum_{j=0}^{m}\binom{n}{k}H_{k,\lambda}(x)H_{\lambda}(m,j)\langle j+m\lambda\rangle_{n-k,\lambda} x^{j}e^{x}. \nonumber
\end{align}
Therefore, by \eqref{39}, we obtain the following theorem. 
\begin{theorem}
For $m,n\ge 0$, we have 
\begin{displaymath}
H_{n+m,\lambda}(x)=\sum_{k=0}^{n}\sum_{j=0}^{m}\binom{n}{k}H_{k,\lambda}(x)H_{\lambda}(m,j)\langle j+m\lambda\rangle_{n-k,\lambda} x^{j}.
\end{displaymath}
\end{theorem}
For sufficiently nice functions $f$, the Fourier inversion theorem says:
\begin{equation}
f(x)=\frac{1}{2\pi}\int_{-\infty}^{\infty}\hat{f}(t)e^{ixt}dt, \label{40}	
\end{equation}
where $\hat{f}(t)=\int_{-\infty}^{\infty}f(x)e^{-ixt}dx$. \\
From \eqref{38} and \eqref{40}, we note that 
\begin{align}
\langle xD\rangle_{n,\lambda}f(x)&=\frac{1}{2\pi}\int_{-\infty}^{\infty}\hat{f}(t)\langle xD\rangle_{n,\lambda}e^{ixt}dt \label{41} \\
&=\frac{1}{2\pi}\int_{-\infty}^{\infty}\hat{f}(t)e^{ixt}H_{n,\lambda}(ixt)dt \nonumber\\
&=\sum_{k=0}^{n}H_{\lambda}(n,k)x^{k}\frac{1}{2\pi}\int_{-\infty}^{\infty}(it)^{k}\hat{f}(t)e^{ixt}dt\nonumber\\
&=\sum_{k=0}^{n}H_{\lambda}(n,k)x^{k}\bigg(\frac{d}{dx}\bigg)^{k}\bigg[\frac{1}{2\pi}\int_{-\infty}^{\infty}\hat{f}(t)e^{ixt}dt\bigg] \nonumber\\
&=\sum_{k=0}^{n}H_{\lambda}(n,k)x^{k}D^{k}f(x).\nonumber
\end{align}
By \eqref{41}, we get the following normal ordering result:
\begin{equation*}
\langle xD\rangle_{n,\lambda}=\sum_{k=0}^{n}H_{\lambda}(n,k)x^{k}D^{k},\quad (n\ge 0).
\end{equation*} \par
For $r\in\mathbb{N}$, we define the {\it{heterogeneous $r$-Stirling numbers of the second kind}} by 
\begin{equation}
\langle x+r\rangle_{n,\lambda}=\sum_{k=0}^{n}H_{\lambda}^{(r)}(n+r,k+r)(x)_{k},\quad (n\ge 0).\label{42}
\end{equation}
Note that 
\begin{align}
e_{\lambda}^{-(x+r)}(-t)&=\sum_{n=0}^{\infty}\langle x+r\rangle_{n,\lambda}\frac{t^{n}}{n!}=\sum_{n=0}^{\infty}\sum_{k=0}^{n}H_{\lambda}^{(r)}(n+r,k+r)(x)_{k}\frac{t^{n}}{n!} \label{43}\\
&=\sum_{k=0}^{\infty}\sum_{n=k}^{\infty}H_{\lambda}^{(r)}(n+r,k+r)\frac{t^{n}}{n!}(x)_{k}.\nonumber
\end{align}
On the other hand, by binomial expansion, we get 
\begin{align}
e_{\lambda}^{-(x+r)}(-t)&=e_{\lambda}^{(-r)}(-t)e_{\lambda}^{-x}(-t)=e_{\lambda}^{-r}(-t)\Big(e_{\lambda}^{-1}(-t)-1+1\Big)^{x} \label{44}\\
&=\sum_{k=0}^{\infty}\frac{1}{k!}\Big(e_{\lambda}^{-1}(-t)-1\Big)^{k}e_{\lambda}^{-r}(-t)(x)_{k}.\nonumber
\end{align}
By \eqref{43} and \eqref{44}, we get 
\begin{equation}
\frac{1}{k!}\Big(e_{\lambda}^{-1}(-t)-1\Big)^{k}e_{\lambda}^{-r}(-t)=\sum_{n=k}^{\infty}H_{\lambda}^{(r)}(n+r,k+r)\frac{t^{n}}{n!}. \label{45}
\end{equation}
By \eqref{42}, we easily get 
\begin{equation}
\begin{aligned}
H_{\lambda}^{(r)}(n+1+r,k+r)&=H_{\lambda}^{(r)}(n+r,k+r-1)\\
&\quad +(k+r+n\lambda)H_{\lambda}^{(r)}(n+r,k+r),
\end{aligned}\label{46}
\end{equation}
where $n,k\in\mathbb{N}$. \par 
We define the {\it{heterogeneous $r$-Bell polynomials}} as
\begin{equation}
H_{n,\lambda}^{(r)}(x)=\sum_{k=0}^{n}x^{k}H_{\lambda}^{(r)}(n+r,k+r),\quad (n\ge 0). \label{47}
\end{equation}
From \eqref{45} and \eqref{47}, we note that 
\begin{equation}
e^{x(e_{\lambda}^{-1}(-t)-1)}e_{\lambda}^{-r}(-t)=\sum_{n=0}^{\infty}H_{n,\lambda}^{(r)}(x)\frac{t^{n}}{n!}. \label{48}
\end{equation}
When $x=1,\ H_{n,\lambda}^{(r)}=H_{n,\lambda}^{(r)}(1)$ are called {\it{heterogeneous $r$-Bell numbers}}. \\
Finally, we observe that 
\begin{align}
\sum_{n=0}^{\infty}H_{n+1,\lambda}(x)\frac{t^{n}}{n!}&=\frac{d}{dt}\bigg(\sum_{n=0}^{\infty}\sum_{k=0}^{n}H_{\lambda}(n,k)x^{k}\frac{t^{n}}{n!}\bigg) \label{49}\\
&=\sum_{n=0}^{\infty}\sum_{k=1}^{n+1}x^{k}H_{\lambda}(n+1,k)\frac{t^{n}}{n!}\nonumber\\
&=\sum_{n=0}^{\infty}x\sum_{k=0}^{n}x^{k}H_{\lambda}(n+1,k+1)\frac{t^{n}}{n!}.\nonumber
\end{align}
Thus, from \eqref{49}, we have 
\begin{align}
H_{n+1,\lambda}(x)&=x\sum_{k=0}^{n}x^{k}H_{\lambda}(n+1,k+1) \label{50}\\
&=\sum_{k=1}^{n+1}H_{\lambda}(n+1,k)x^{k}\quad (n\ge 0). \nonumber
\end{align}
\section{Conclusion}
In conclusion, this paper successfully introduced and explored the heterogeneous Stirling numbers of the second and first kinds, $H_{\lambda}(n,k)$ and $G_{\lambda}(n,k)$, respectively, along with a polynomial extension of $H_{\lambda}(n,k)$, the heterogeneous Bell polynomials $H_{n,\lambda}(x)$, and their $r$-analogues. These new number sequences and polynomials provide a unifying framework that generalizes both the classical Stirling numbers and the Lah numbers, seamlessly transitioning between these well-established combinatorial objects as the parameter $\lambda$ varies. \par
We derived fundamental properties, including generating functions, explicit expressions, and recurrence relations, for each of these newly defined sequences. The establishment of a Dobinski-like formula for the heterogeneous Bell polynomials and the discovery of combinatorial identities further enriched our understanding of their combinatorial significance. Moreover, the extension to the $r$-case demonstrated the robustness and versatility of these definitions, opening avenues for further exploration in related areas. \par
Future research could focus on exploring deeper combinatorial interpretations, connections to other special functions, and applications in various scientific and engineering disciplines. The heterogeneous Stirling and Bell numbers, with their rich mathematical structure, hold promise for further discoveries and applications.

\end{document}